\documentclass[12pt]{article}
\usepackage[top=3.3cm,bottom=3.2cm,left=2.8cm,right=2.8cm]{geometry}
\usepackage[utf8]{inputenc}
\usepackage{amssymb,amsmath,amsfonts,eucal,mathrsfs,amsthm,xcolor} 
\usepackage{enumitem}
\usepackage{amsrefs}
\usepackage{stmaryrd}
\usepackage{hyperref}
\hypersetup{
    colorlinks=true,
    filecolor=magenta,      
    linkcolor={blue!50!black},
    citecolor={red!50!black},
    urlcolor={blue!80!black}
}
\usepackage{dutchcal}

\newtheorem*{theorem*}{Theorem}
\newtheorem{theorem}{Theorem}
\newtheorem{proposition}[theorem]{Proposition}

\theoremstyle{definition}

\newcommand{\R}{\mathbb{R}}

\newcommand{\Q}{\mathbb{Q}}
\newcommand{\Sf}{\mathbb{S}}

\newcommand{\Hy}{\mathbb{H}}

\newcommand{\hess}{\mbox{Hess\,}}

\newcommand{\grad}{\mbox{grad\,}}

\newcommand{\ind}{{\rm Index}\, }

\def\span{{\rm{span}}}
\def\bea{\begin{eqnarray*} }
\def\eea{\end{eqnarray*} }

\def\beq{\begin{equation}}
\def\Z{\mathord{\mathbb Z}}

\def\B{\mathcal{B}}
\def\<{{\langle}}
\def\>{{\rangle}}

\def\n{\nabla}

\def\a{\alpha}

\def\be{\begin{equation} }
\def\ee{\end{equation} }

\begin{document}

\title{On the first eigenvalue of the Hodge Laplacian \\ of submanifolds}
\author{Christos-Raent Onti}
\date{}
\maketitle

\renewcommand{\thefootnote}{\fnsymbol{footnote}} 
\renewcommand{\thefootnote}{\arabic{footnote}} 

\renewcommand{\thefootnote}{\fnsymbol{footnote}} 
\renewcommand{\thefootnote}{\arabic{footnote}}

\begin{abstract}
We prove that equality in a sharp lower bound for the first $p$-eigenvalue 
of the Hodge Laplacian on closed submanifolds in space forms can occur 
only on topological spheres, assuming positivity.
\end{abstract}

\section{Introduction}

Let $M^n$ be a closed, connected and oriented Riemannian manifold of 
dimension $n$. For each integer $1\leq p\leq n-1$, the Hodge-Laplace 
operator (or the Hodge Laplacian) acting on $p$-forms is defined by 
$$
\Delta=d\delta+\delta d: \Omega^p(M^n)\to \Omega^p(M^n),
$$
where $d$ and $\delta$ are the differential and the co-differential operators, respectively. 
It is well known that the spectrum of the Hodge-Laplace operator is discrete and non-negative, 
and that its kernel is isomorphic to the $p$-th de Rham cohomology group $H^p(M^n;\R)$.
If $\lambda_{1,p}(M^n)$ denotes its lowest eigenvalue, then 
$$
 \lambda_{1,p}(M^n)=\inf_{\omega\in \Omega^p(M^n)\setminus\{0\}}\frac{\int_{M}\left(\|d\omega\|^2+\|\delta\omega\|^2\right) dM}{\int_M \|\omega\|^2  dM}.
$$
Since the above is invariant by the Poincar\'{e} duality induced by the Hodge $*$-operator, we have $\lambda_{1,p}(M^n)=\lambda_{1,n-p}(M^n)$ 
and thus we may assume that $p\leq n/2.$ Moreover, it is clear that if $\lambda_{1,p}(M^n)>0$, then $H^p(M^n;\R)=H^{n-p}(M^n;\R)=0$.

The Hodge Laplacian satisfies for every $p$-form $\omega\in\Omega^p(M^n)$ the Bochner-Weitzenb\"ock formula
\be\label{bwform}
\Delta\omega=\n^*\n\omega+\B^{[p]}\omega,
\ee
where $\n^*\n$ is 
the connection Laplacian and 
$
\B^{[p]}\colon \Omega^p(M^n)\to \Omega^p(M^n)
$
is a certain symmetric endomorphism on the bundle of $p$-forms, called the {\it Bochner-Weitzenb\"ock operator}.
Therefore, \eqref{bwform} implies that lower bounds on the Bochner-Weitzenb\"ock operator lead naturally to lower bounds on the Hodge-Laplace operator. 
In particular, from \cite[Proposition 3]{Savo14} we get that
\begin{equation}\label{sala14}
\text{if}\; \B^{[p]}\geq p(n-p)\Lambda\; \text{for some} \; \Lambda>0,\;  \text{then}\; \lambda_{1,p}(M^n)\geq p(n-p+1)\Lambda.
\end{equation}

Let $f\colon M^n\to \tilde M^{n+m}, n\geq 3,$ be an isometric immersion into a Riemannian manifold $\tilde M^{n+m}$ of dimension $n+m$.
The second fundamental form $\a_f$ is viewed as a section of the vector bundle $\mathrm{Hom}(TM\times TM,N_f M)$, where 
$N_f M$ is the normal bundle. For each unit normal vector field $\xi\in \Gamma(N_fM)$, the associated 
shape operator $A_{\xi}$ is given by 
$$
\<A_{\xi} X,Y\>=\langle \alpha_f(X,Y),\xi\rangle,\: \: X,Y\in TM.
$$
Recall that the traceless part of the second fundamental form is given by $\Phi=\a_f-\<\cdot,\cdot\>\mathcal H$, where $\mathcal H$ denotes the
\textit{mean curvature vector field} given by $\mathcal H=(\mathrm{tr}\, \alpha_f)/n$, where $\mathrm{tr}$ 
means taking the trace. Finally, by $H$ we denote the length of the \textit{mean curvature}, that is, $H=\| \mathcal H \|$.
In \cite[Proposition 16]{ov22} we proved with Vlachos that

\begin{proposition}\label{bprop}
If the curvature operator of $\tilde M^{n+m}$ is 
bounded from below by a constant $c$, then the Bochner 
operator of $M^n$, for any $1\leq p\leq \lfloor n/2\rfloor$, satisfies pointwise the inequality
\be\label{bineq}
\mathop{\min_{\omega\in \Omega^p(M^n)}}_{\|\omega\|=1}\<\B^{[p]}\omega,\omega\>
\geq \frac{p(n-p)}{n} \big(n(H^2+c) -\frac{n(n-2p)}{\sqrt{np(n-p)}}\, H\|\Phi\|-\|\Phi\|^2\big).
\ee
If equality holds in \eqref{bineq} at a point $x\in M^n$, then the following hold:
\begin{enumerate}[topsep=2pt,itemsep=2pt,partopsep=2ex,parsep=0.5ex,leftmargin=*, label=(\roman*), align=left, labelsep=-0em]
\item The shape operator $A_\xi(x)$ has at most two distinct 
eigenvalues with multiplicities $p$ and $n-p$ for every unit vector $\xi\in N_fM(x)$. If in addition 
$p<n/2$ and the eigenvalue of multiplicity $n-p$ vanishes, then $A_\xi(x)=0$.
\item If $H(x)\neq 0$ and $p<n/2$, then $\mathrm{Im}\, \a(x)=\span\left\{\mathcal H(x)\right\}$.
\end{enumerate}
\end{proposition}

\noindent Therefore, if
$$
\kappa_p:=\min_{x\in M^n}
\left\{(H^2+c) -\frac{n-2p}{\sqrt{np(n-p)}}\, H\|\Phi\|-\frac{1}{n}\|\Phi\|^2\right\}
$$
for some $1\leq p\leq \lfloor n/2\rfloor$, 
then it follows from \eqref{sala14} and \eqref{bineq} that
\be\label{feigen}
\lambda_{1,p}(M^n)\geq p(n-p+1)\kappa_p.
\ee

Inequality \eqref{feigen} was first proved by Savo for hypersurfaces \cite[Theorem 7]{Savo14}, and 
subsequently extended by Cui and Sun to submanifolds of arbitrary codimension \cite[Theorem 1.1]{quisun19}. 
They also showed that the inequality is sharp by providing trivial examples attaining equality. However, 
no characterization was given of the submanifolds for which equality holds.
The aim of this note is to shed light on the case of equality in \eqref{feigen} assuming $\lambda_{1,p}(M^n)>0$, 
when $\tilde M^{n+m}=\Q_c^{n+m}$, where $\Q_c^{n+m}$ denotes the complete simply 
connected space form of constant sectional curvature $c$. In fact, we prove that in this case 
equality occurs only on topological spheres.
For simplicity we assume that 
$c\in\{0,\pm1\}$. Thus  $\Q_c^{n+m}$ is the Euclidean space $\R^{n+m} (c=0),$ the unit 
sphere $\Sf^{n+m} (c=1),$ or the hyperbolic space $\Hy^{n+m}\; (c=-1)$. 

\begin{theorem*}
Let $f\colon M^n\to  \Q_c^{n+m}, n\geq 4$, be an isometric immersion of a closed, connected and oriented Riemannian manifold. 
If for some $1\leq p\leq \lfloor n/2\rfloor $ equality holds in \eqref{feigen} with $\lambda_{1,p}(M^n)>0$, then $M^n$ is homeomorphic to the 
sphere $\Sf^n$.
\end{theorem*}

\section{Proof of the Theorem}
The idea of the proof is to show that $M^n$ is a simply connected homology sphere over the integers
and the proof will follow by the generalized Poincar\'e conjecture (Smale $n\geq 5$, Freedman $n=4$).

Assume that for some $1\leq p\leq \lfloor n/2\rfloor$ equality holds in \eqref{feigen} with $\lambda_{1,p}(M^n)>0$. 
Then Proposition \ref{bprop}  implies that the shape operator $A_\xi(x)$ at each point $x$ will have at most two distinct 
eigenvalues of multiplicities $p$ and $n-p$ for every unit vector $\xi\in N_fM(x)$. 
We claim that there exists a Morse function on $M^n$ such that 
the index at each critical point is $0,p,n-p$ or $n$. To this end, we distinguish the following two cases:

\medskip

\noindent {\underline {\sc Case $c\in\{0,1\}$}}: Let $u\in\R^{n+m+c}$ be a
vector such that the height function
$$
\varphi\colon M^n\to \R, \; \varphi(x)=\<f_c(x),u\>
$$
is a Morse function, where 
$$
f_c=\begin{cases}
f, & \text{if} \;  c=0, \\[2mm]
j\circ f, & \text{if} \; c=1, \;\text{and} \; j\colon \Sf^{n+m}\to \R^{n+m+1} \; \text{denotes the standard inclusion}.
\end{cases}
$$
A direct computation gives that at a critical point $x_0$ of $\varphi$ we have
$$
u\in N_{f_c}M(x_0)\; \text{and}\; \hess \varphi (X,Y)=\<\alpha_{f_c}(X,Y),u\>,\; \text{for all}\; X,Y\in T_{x_0}M.
$$
Obviously, the second fundamental form of $f_c$ has at most two distinct principal curvatures of 
multiplicities $p$ and $n-p$ in every normal direction and the claim follows in this case.

\medskip

\noindent  {\underline {\sc Case $c=-1$}}: 
We consider the 
function $$\varphi\colon \Hy^{n+m}\to \R,\; \varphi(x)=\frac{1}{2}r^2(x),$$ 
where $r(x)$ denotes the distance function issuing from some suitable choice of 
point $o\in\Hy^{n+m}$ to $x\in\Hy^{n+m}$. It is a standard fact that $\varphi$ is smooth inside the cut locus of $o$. 
Let $\gamma(t)$ be a unit speed geodesic with $\gamma(0)=o$. Then, 
we have $\gamma^\prime(t)=\grad r(\gamma(t))$.
For $X,Y\in \Gamma(T\Hy)$ a direct computation gives 
$$
\<X,\grad \varphi\>=r\<X,\grad r\>
\;\text{and}\;
\hess \varphi(X,Y)=\<X,\grad r\>\<Y,\grad r\> +r\hess r(X,Y).
$$
Consider $\tilde \varphi=\varphi\circ f\colon M^n\to (0,+\infty)$ and choose $o\in\Hy^{n+m}$ such that $\tilde \varphi$ is 
a Morse function on $M^n$ (notice that this is always possible). 
At a critical point $x_0\in M^n$, we have $\grad r(f(x_0))\perp f_*(T_{x_0} M^n)$, that is the (unique unit speed) 
geodesic $\gamma(t)$ will hit $f(M^n)$ orthogonally, and
\be\label{subhes1}
\hess \tilde\varphi(\tilde X,\tilde Y)=r\left(\hess r(f_*\tilde X, f_*\tilde Y)+\<A_{{\rm{grad}}\; r} \tilde X,\tilde Y\> \right), \; \tilde X,\tilde Y\in T_{x_0}M^n.
\ee
Let $\gamma(\ell)=f(x_0)$ and consider a Jacobi field $J(t)$ along $\gamma(t)$ with $J(0)=0$. It follows that
\be\label{rhes1}
\hess r(J(\ell),J(\ell))=\<J^\prime(\ell),J(\ell)\>=\frac{1}{2}\frac{d}{dt}\|J(t)\|^2|_{t=\ell}.
\ee
Recall that the Jacobi fields $J(t)$ in $\Hy^{n+m}$ with $J^\prime(0)\perp \gamma^\prime(0)$ are given by
$$
J(t)=\sinh t\cdot w(t),
$$
where $w(t)$ is a parallel vector field along $\gamma(t)$ with $J'(0)=w$ and $\|w\|=1$.
Observe that $A_{{\rm{grad}}\; r}(x_0)$ has at most two distinct eigenvalues, say $\lambda$ and $\mu$ with 
multiplicities $p$ and $n-p$, respectively.
Consider an orthonormal basis $\{e_1,\dots,e_n\}$ of $T_{x_0} M^n$ such that 
$$
A_{{\rm{grad}}\; r}(e_i)=\lambda e_i,\;\; 1\leq i\leq p,
\;\;\text{and}\;\; 
A_{{\rm{grad}}\; r}(e_i)=\mu e_i,\;\; p+1\leq i\leq n.
$$
Let $w_i(t)$ such that 
$$
w_i(\ell)=f_*(e_i),\; 1\leq i\leq n,
$$ 
with corresponding Jacobi fields 
$$
J_i(t)=\sinh t\cdot w_i(t),\; 1\leq i\leq n.
$$
Therefore, from \eqref{rhes1} we obtain 
$$
\hess r(J_i(\ell),J_i(\ell))=\frac{1}{2}\frac{d}{dt}\|J_i(t)\|^2|_{t=\ell}=
\displaystyle \frac{1}{2}\sinh(2\ell),\; \text{for all}\; 1\leq i \leq n.
$$
Hence \eqref{subhes1} gives 
$$
\hess \tilde\varphi(e_i,e_i)=
\begin{cases}
\displaystyle \ell(\coth \ell+ \lambda), & \;\; \text{for}\;\; 1\leq i\leq p,\\[4mm]
\displaystyle \ell(\coth \ell+ \mu), & \;\; \text{for}\;\; p+1\leq i\leq n,
\end{cases}
$$
and therefore it is now clear that $\ind\tilde\varphi(x_0)\in\{0,p,n-p,n\}$. This completes the proof of 
the claim.

Therefore, it follows from standard Morse theory 
(cf. \cite[Th. 3.5]{Milnor63} or \cite[Th. 4.10]{CE75}) that 
$M^n$ has the homotopy type of a CW-complex with cells only in dimensions $0,p,n-p$ or $n$.
Therefore 
\be\label{hgroup}
H_i(M^n;\Z)=0,\; \text{for all}\; i\neq 0,p,n-p, n.
\ee 
Next, we claim that also
\be\label{pgroup}
H_p(M^n;\Z)=H_{n-p}(M^n;\Z)=0.
\ee
Indeed, our hypothesis 
implies
$$
H^p(M^n;\R)=H^{n-p}(M^n;\R)=0.
$$
Hence
\be\label{hgroup2}
H_p(M^n;\Z)={\rm Tor}(H_p(M^n;\Z))
\;\; \text{and}\;\; 
H_{n-p}(M^n;\Z)={\rm Tor}(H_{n-p}(M^n;\Z)).
\ee
By the Poincar\'{e} duality, the universal coefficient theorem and \eqref{hgroup}, we have 
$$
{\rm Tor}(H_p(M^n;\Z))\cong {\rm Tor}(H_{n-p-1}(M^n;\Z))=0
$$
and
$$
{\rm Tor}(H_{n-p}(M^n;\Z))\cong {\rm Tor}(H_{p-1}(M^n;\Z))=0, 
$$ 
where $\cong$ denotes the isomorphism of the corresponding groups. This, in combination with \eqref{hgroup2},
proves \eqref{pgroup}. Hence $M^n$ 
is a homology sphere over the integers.

Finally, we show that $M^n$ is simply connected. If $p\neq 1$ this follows directly from 
 \cite[Proposition 4.5.7, p. 90]{AD}, as in this case, $\varphi$ has no critical points of index one.
If $p=1$, then since there are no $2$-cells, it follows by the cellular approximation theorem 
that the inclusion of the $1$-skeleton $\mathrm{X}^{(1)}\hookrightarrow M^n$ 
induces isomorphism between the fundamental groups. Therefore, 
$\pi_1(M^n)$ is a free group on $\beta_1(M^n;\Z)=0$ elements, and thus
$M^n$ is simply connected.

Therefore, $M^n$ is a homotopy sphere and by the generalized
Poincar\'e conjecture (Smale $n\geq 5$, Freedman $n=4$), 
$M^n$ is homeomorphic to $\Sf^n$, which concludes the proof of the theorem.

\vspace{10mm}

\noindent Christos-Raent Onti\\
Department of Mathematics and Statistics\\
University of Cyprus\\
1678, Nicosia -- Cyprus\\
e-mail: onti.christos-raent@ucy.ac.cy
\bigskip

\end{document}